\documentclass{article} 
\usepackage{a4wide}
 
\newcommand{\np}{\relax}
\newcommand{\mx}{{\max}}

\usepackage{eepic}
\usepackage{amssymb}
\usepackage{amsthm}
\usepackage{mathrsfs}
\usepackage{graphicx}
\usepackage{color}

\def\on{{\upharpoonright}}

\def\itm#1 {\item[{(#1)}]}

\newenvironment{ilist}{\begin{list}{(\arabic{enumi})}{\usecounter{enumi}
\leftmargin=2cm\labelsep=3mm\labelwidth=3cm}}{\end{list}}

\newtheorem{thm}{Theorem} [section]

\newtheorem{Theorem}[thm]{Theorem}
\newtheorem{lemma}[thm]{Lemma}
\newtheorem{Lemma}[thm]{Lemma}

\newtheorem{conclusion}[thm]{Conclusion}

\newtheorem{Corollary}[thm]{Corollary}

\theoremstyle{remark}

\newtheorem{fact}[thm]{Fact}
\newtheorem{Fact}[thm]{Fact}
\newtheorem{example}[thm]{Example} 

\newtheorem{examples}[thm]{Examples}

\newtheorem{dfact}[thm]{Fact and Definition}

\newtheorem{remark}[thm]{Remark}
\newtheorem{Remark}[thm]{Remark}
     
\newtheorem{defn}[thm]{Definition}
\newtheorem{definition}[thm]{Definition}     
\newtheorem{Definition}[thm]{Definition}     
\theoremstyle{definition} 

\newtheorem{notation}[thm]{Notation}
\newtheorem{Notation}[thm]{Notation}

\newtheorem{motivation}[thm]{Motivation}

\newcommand{\cl}[1]{\langle #1 \rangle}
\newcommand{\clmax}[1]{\cl{#1}_\mx}

\newcommand{\rem}[1]{\setbox0\hbox{#1}\ifdim\wd0<0.7\hsize 
\par\hbox to \hsize{\hfil\vrule\hskip2pt\vrule\ \box0\ \vrule\hskip2pt\vrule}
 \else
\par\hbox to \hsize{\hfil \vrule\hskip2pt\vrule\ 
\vbox{\parindent=0cm\hsize=0.7\hsize \strut\sl #1\strut}\ \vrule\hskip 2pt\vrule}\par\fi}

\def\xx#1^#2 {#1^{(#2)}}  
 \def\yy#1^#2{#1^{[#2]}}  
 \def\yyy#1^#2<{#1^{[{\rm root}>#2]}}  

\let\zz\yy
\def\uu#1_#2{#1_{\langle#2\rangle}}
 \newcommand{\intr}{{\bf int}}

 \newcommand{\leaf}{{\bf ext}}   
 \newcommand{\ext}{\leaf}
\newcommand{\llim}{{\bf lim}}

  \newcommand{\wurzel}{{\bf root}} 
 \newcommand{\Wurzel}{{\bf Root}} 
\renewcommand{\succ}{{\rm Succ}}
\newcommand{\Succ}{{\rm Succ}}

\newcommand{\C}{{\mathscr C}}
\newcommand{\D}{{\mathscr D}}

\renewcommand{\O}{{\mathscr O}}

\begin{document}

\parindent=0cm

\renewcommand{\qedsymbol}{//}

\newcommand{\CAP}{{\cap}}
\renewcommand{\baselinestretch}{1.4}

\newcommand{\drop}{_{\rm drop}}
\newcommand{\thin}{_{\rm thin}}
\newcommand{\short}{_{\rm short}}
\newcommand{\thinl}{_{\rm thin}}

\newcommand{\thinshort}{_{\rm thin/short}}
\newcommand{\glue}{_{\rm glue}}
\newcommand{\DROP}{{\it DROP}}
\newcommand{\THIN}{{\it THIN}}
\newcommand{\GLUE}{{\it GLUE}}
\newcommand{\enull}{{\eta_0}}
\newcommand{\es}{_{\ext(S)}}

\newcommand{\rda}{{r(\delta\cdot\xi)}}
\newcommand{\rdai}{{r(\delta\cdot\xi+i)}}
\newcommand{\dai}{{\delta\cdot\xi+i}}
\newcommand{\rdaa}[1]{{r(\delta\cdot\xi+#1)}}

\def\oo#1 {\O^{(#1)}}
\def\cid{{\C}_{\rm id}}

\title{All creatures great and small} 
\author{Martin Goldstern \and Saharon Shelah}

\date{2006-02-xx, 2007-06-08} 

\pagestyle{headings}

\maketitle
\begin{abstract}
Let $\kappa$ be an uncountable regular cardinal.
Assuming $2^\kappa=\kappa^+$, we show that the
clone lattice on a set of size $\kappa$ is not dually atomic.
\end{abstract}


\setcounter{section}{-1} 

\section{Introduction }  
  
 A clone $\C$ on a set $X$ is a set of finitary operations $f:X^n\to X$ 
 which contains all the projections and is closed under 
 composition. 
 (Alternatively, $\C$ is a clone {if} $\C$ is the set of term functions of 
 some universal  algebra over~$X$.)  
  
  The family of all clones forms a complete algebraic lattice $Cl(X)$
 with greatest element $\O = \bigcup_{n=1}^\infty X^{X^n} $, where $
 X^{X^n} $ is the set of all $n$-ary operations on~$X$.  (In this
 paper, the underlying set $X$ will be a fixed uncountable set.)

   The coatoms of this lattice $Cl(X) $ are called 
 ``precomplete clones'' or ``maximal clones'' on~$X$.  
  
The classical reference  for older results
 about clones  is \cite{PK:1979}. 

 For singleton sets $X$ the lattice $Cl(X)$ is trivial; for $|X|=2$
 the lattice $Cl(X)$ is countable, and well understood (``Post's
 lattice'').   For $|X|\ge 3$,  $Cl(X)$ has uncountably many elements.  
 Many results for clones on finite sets 
  can be found in  
 \cite{Szendrei:1986} and the recent \cite{Lau:2006}.
   In particular, there is an explicit
 description of all  
 (finitely many) precomplete clones on a given finite set 
  (\cite{Rosenberg:1970}, see also   
   \cite{Quackenbush:1971} and \cite{Buevich:1996});
 this description also includes a decision procedure
 for the membership problem for each of these clones. 
It is also known 
 that every clone $\C\not=\O$ is contained in a precomplete 
 clone, that is: the clone lattice  $Cl(X)$  on any finite set  
 $X$ is {\em dually atomic}. 
  (This gives an explicit criterion for deciding whether a given 
 set of functions generates all of~$\O$: just check {if} it is contained 
 in one of the precomplete clones.)  
  
 Fewer results are known about the lattice of clones on an infinite 
 set, and they are often negative or ``nonstructure'' results:
   \cite{Rosenberg:1976} showed that 
 there are always $2^{2^\kappa}$ precomplete clones on a set of 
 infinite cardinality~$\kappa$ (see also \cite{GoSh:737}). 
  
  \cite {Rosenberg+Schweigert:1982} 
 investigated  
 ``local'' clones on infinite sets (clones that are closed sets  in the 
 product topology).  It  is easy to see that the lattice of local 
 clones is far from being dually atomic. 

Already Gavrilov in \cite[page 22/23]{Gavrilov:1959} asked whether the lattice
 of all clones on a {\em countable} set is also dually atomic, since a
 positive answer would be an important component for a completeness
 criterion, as remarked above.  The same question for {\em all} infinite sets 
is  listed as  problem P8 in \cite[page 91]{PK:1979}.

In \cite{GoSh:808} we showed that that (assuming the continuum
hypothesis CH) the answer is negative for countable sets.  We are now
able to extend  
the construction from
\cite{GoSh:808}  to work on all regular uncountable cardinals
as long as they satisfy the corresponding version of CH.   The
question whether such a theorem is provable in ZFC alone remains
open. 

   We will
write CH$_\lambda$  
for the statement $2^\lambda = \lambda^+$, or equivalently,
\begin{quote}
\begin{enumerate}
\item[CH$_\lambda$:\ ]
If $|X| = \lambda$, then every subset of~${\mathscr P}(X)$ (the power set of
$X$)\\  either has cardinality $\le \lambda$, or is equinumerous with
 $ {\mathscr P}(X)$
\end{enumerate}
\end{quote}

We will show here the following for every uncountable regular cardinal $\lambda$: 
 \begin{Theorem}
  \label{thm} 
Assume that $X$ is a set of size $\lambda$, and that CH$_\lambda$ holds. 
\\
Then 
 the lattice of clones on~$X$ 
set is {\bf not dually atomic}, i.e., there is a clone 
 $\C\not=\O$ which is not contained in any precomplete clone.  
 \end{Theorem}

 The clone $\C_U$ that we construct has the additional feature that we can 
 give a good description of the interval~$[\C_U,\O]$.


The method behind our proof is ``forcing with large creatures'', a
new method which is rooted in ``forcing with normed creatures''
(\cite{Sh:207}, \cite{RoSh:470}). 
However, for the purposes of 
 this paper the connection with forcing machinery is sufficiently 
 shallow to allow us to be self-contained.  
  
 In particular, no knowledge  
 of set theory 
  is required for our theorem, except for a basic understanding of
  ordinals, well-founded relations and transfinite induction.


\begin{Remark}
The reader familiar with \cite{GoSh:808} may appreciate the following
list of differences/modifications: 
\begin{enumerate}
\item In \cite{GoSh:808}, our ``largeness property'' was connected
  with cardinalities of finite sets going to infinity, and we could
  show several partition theorems of the form:  if the norms of a
  sequence of creatures $(S_n)$ goes to infinity, we can find a
  subsequence $(T_n)$ of ``nice'' creatures (e.g., homogeneous for some
  colouring function) such that their norm still goes to infinity. 
\\
 This point has become easier now; rather than sets ``large in
  cardinality'', our large sets are now sets in certain
  ultrafilters. 
\item In \cite{GoSh:808} we had ``unary'' and ``binary'' partition
  theorems guaranteeing that we can thin out creatures to creatures that are homogeneous with respect to certain coloring functions.  
In the current paper we only have a unary partition
  theorem (see \ref{unary}).
 This means that our notions of ``$f$-weak'' and ``$f$-strong''
  are somewhat weaker than the notions in
 \cite{GoSh:808}, which in turn means that we know somewhat less 
about the structure of the clone interval we construct.   In particular, 
instead of showing that this interval is linearly ordered, we can only 
show that there is a linearly ordered cofinal set. 

\item In  \cite{GoSh:808}, our construction took $\omega_1$ steps, so
 in each intermediate step we only had to consider the countably many
 steps taken before. In particular, the ${\sigma}$-closure of our 
 set of creatures in  \cite{GoSh:808} was easily proved via a
 ``diagonal'' construction.   In the current paper we again have a
simple diagonal construction (\ref{fusion.lambda}) 
to find a lower bound of a decreasing chain of creatures of
 length $\lambda$, but we also have to deal with shorter infinite
 sequences in \ref{fusion.delta}, which necessitates a more complicated
 setup. 
\item For any $f:\lambda \to \lambda $ let $\bar f: \lambda \to
  \lambda $ be defined as $\bar f(x) = \sup\{f(x): x \le y \}$. \\
If $\lambda = \omega$, then we have $f\in \C$ iff $\bar f\in \C$ for
  all (relevant) clones $\C$, so in  \cite{GoSh:808} we could wlog assume that all unary
 functions that we considered were monotone.   But for $\lambda > \omega$ we cannot
  assume that any more. 
\item We introduce ``coordinates'' for elements of creatures.  This will obviate the notational difficulties we had in \cite[3.10]{GoSh:808} (involving the 
possible ``recycling'' of deleted notes). 
\item Another notational change: rather than
defining a linear order of equivalence classes of fronts as
in \cite[5.2]{GoSh:808}, we will work directly with the induced order
on the functions in~$\O$. 
\end{enumerate}
\end{Remark}

\np

\section{Preliminaries}

Our base set will be a fixed uncountable regular cardinal $\lambda$, equipped with the usual order. We 
are interested in operations on~$\lambda$, i.e., elements of 
$\O = \bigcup_{k=1,2,\ldots} \lambda^{\lambda^k}$, and in subsets of~$\O$.

\begin{definition}
We write $\C_\mx $ for the set  of all functions $f$ which satisfy 
$f(x_1,\ldots, x_k) \le \max(x_1,\ldots, x_k)$ for all $x_1,\ldots, x_k\in \lambda$.

For each set $\D \subseteq \O$ we write $\cl{\D}$ for the clone generated
by $\D$.
We will write $\clmax{\D}$ for $\cl{\C_\mx\cup \D}$.
\end{definition}

\begin{fact}\label{c.max}
\begin{enumerate}
\item 
 $\C_\mx$ is a clone.
\item 
 Any clone containing $\C_\mx$ is downward closed (in the sense of the
pointwise partial order on each of the sets $\lambda^{\lambda^n}$). 
\item 
Assume that  $\C\supseteq \C_\mx$ is a clone, and assume that 
 $f_1, \ldots, f_k$ are functions of the same arity. 
Then 
 $\cl{\C\cup \{f_1,\ldots, f_k\}} = \cl{\C\cup \{\max(f_1,\ldots, f_k\})}$.
\\ (Here, $\max$ is the pointwise maximum function.)
\end{enumerate}
\end{fact}

\begin{proof} 
(1) is trivial, and (2) is easy (see \cite{GoSh:808}):  If $g\in \C$, and 
$f$ is $k$-ary, 
$f(\vec x) \le g(\vec x)$ for all $\vec x$, then we can find a $k+1$-ary 
function $F\in \C_\mx$ with $f(\vec x) = F(\vec x, g(\vec x))$ for all $\vec x$.

In (3), the inclusion
    $\subseteq$ follows from the downward closure of
 $\cl{\C\cup \{f_1,\ldots, f_k\}}$ and (2), and the inclusion $\supseteq$ follows from the assumption that the $k$-ary maximum function is in $\C$. 
\end{proof}

\np

\subsection{Proof outline}\label{proof.outline}

 \begin{Fact}  Let $(L,{<})$
  be a complete linear order.  $\C\supseteq \C_\mx$ a clone, and let
 $\rho:\O\to L$ be a map into $L$ with properties (a), (b), (c) (where
 we write $f<_\rho g$ for $\rho(f)< \rho(g)$, similarly $f\le_\rho
 g$).  Then (1), (2), (3) hold.
\begin{itemize}
\itm a $ f <_\rho g \ \Rightarrow \ f\in \cl {\C\cup\{g\}}$.
\itm b  $f\in \cl {\C\cup\{g\}}  \ \Rightarrow \ f\le_ \rho g$. 
\itm c $\rho(\max(f, g)) = \max(\rho(g), \ldots, \rho(g))$. 
\end{itemize}
 Then 
 \begin{enumerate}
 \itm 1 For every $d\in   L$ the sets $\D_{<d}:=\{f: \rho(f)<d\}$ and 
 $\D_{\le d}:=\{f: \rho(f)\le d\}$ are clones (unless they are  empty).
 \itm 2  For every clone $\D$ in  $[\C, \O)$ there is some $d\in L$ with 
 $\D_{<d} \subseteq \C \subseteq \D_{\le d}$. 
 \itm 3 If moreover $\rho[\O]$ has no last element, then the interval $[\C,\O]$
 has no coatom. 
 \end{enumerate}
Note that $\forall f\, \forall g\, $(b) is equivalent to $\forall f\,
\forall g\, $(b'), and (a)+(b)+(c)  is equivalent to (a)+(b')+(c'):
\begin{enumerate}
\itm b'  $f<_\rho g  \ \Rightarrow \ g\notin \cl{\C\cup \{f\}}$. 
\itm c' Whenever $f<\rho g$ or $f\sim_\rho g$, then $\max(f,g)\sim_\rho g$.
\end{enumerate}
 \end{Fact}

\begin{proof}
Writing $0$ for $\inf \rho[\O]$, we conclude from (b): 
$$ f\in \C \Rightarrow \rho(f)=0$$

Property (c) implies that the sets $\D_{<e}$ and $\D_{\le e}$
 are closed under the pointwise $\max$ function; if they are nonempty, they
contain $\C$ (and hence also all projections). 

For $e\ge 0$, $k>0$ we show $\cl{f_1,\ldots, f_k} \subseteq \D_{\le e}$ for any
 $f_1,\ldots, f_k\in \D_{\le e}$:
\begin{quote}
Let $f:=\max(f_1,\ldots, f_k)\in \D_{\le e}$.
 So $\cl{\C \cup \{f_1,\ldots, f_k\}} = \cl{\C \cup \{f\}}$. \\
If $h\in \cl{\C\cup \{f\}}$, then (by (b)) $\rho(h)\le \rho(f) \le e$.  So
$\cl{\C\cup \{f\}} \subseteq \D_{\le e}$. 
\end{quote}

Hence $\D_{\le e}$
is a clone. 

The argument for $\D_{< e} $ (with $e>0$) is similar.

Now, given any clone $\D\supseteq \C$, let $d_0:= \sup \{\rho(f):f\in
\D\}$.   We claim $\D_{<d_0} \subseteq \D \subseteq \D_{\le d_0}$: 
\begin{quote}
Clearly $\D \subseteq \D_{\le d_0}$. \\ Let $h\in \D_{< d_0}$,
then $\rho(h)<d_0$, hence there is some $f\in \D$ with
$\rho(h)<\rho(f)$.  So $h\in \cl{\C\cup \{f\}}\subseteq \D$ by (a).
Hence $\D_{<d_0} \subseteq \D$.
\end{quote}

Finally, we see that the map $d\mapsto \D_{\le d}$ is 1-1 from $\rho[\O]$
into $[\C,\O)$, since $\rho(f)=d$ implies $f\in \D_{\le d}\setminus \D_{\le e}$ for $e<d$. Hence $[\C,\O)$ contains a cofinal copy of $\rho[\O]$, thus no 
maximal  element. 
\end{proof}

We will try to find a linear order $L$ and a map $\rho$ that will allow us 
to apply the lemma.  But rather than finding $L$ explicitly, we will
first 
construct relations  $<_\rho$ and $\sim_\rho$
$$(**) \qquad \qquad 
f <_\rho g \Leftrightarrow \rho(f)<\rho(g) \qquad \qquad 
f \sim_\rho  g \Leftrightarrow \rho(f)=\rho(g)$$
on $\O$.  The order $L$ will then appear as the Dedekind completion of the quotient order $\O/{\sim}$. 

We will construct $<$ and $\sim $ in $\lambda^+$ many stages, as unions 
$\bigcup_i {<}_i$ and $\bigcup_i {\sim}_i$.    Each $<_i$ will be a partial order on $\O$, and each $\sim_i$ will be an equivalence relation, but only at the end will we guarantee that any two operations $f$ and $g$ are either $<$-comparable or $\sim$-equivalent. 

The relation $f<_i g$ will say that on a ``large'' set, $f$ grows
faster than $g$. This $i$-th notion of ``large'' will come from a
filter $D_i$ on $\lambda$.  Eventually, the clone $\C$ at the bottom
of our interval will be determined by the filter $\bigcup_i D_i$.

\np
\subsection{Filter clones} 

\begin{definition}\label{h.A}
For any unbounded  $A \subseteq \lambda $, let $h_A$ be the function 
$h_A(x) = \min \{y\in A: y > x\}$. 

For any family $U$ of unbounded subsets of~$\lambda$ let $\C_U$ be 
the clone $\clmax{ h_A:A\in U} $.
\end{definition}

(The function $h_F$ will be defined below in \ref{h.F}.)

\begin{lemma}\label{U.h.A}
  Assume that $U$ is a filter on $\lambda$ containing no
 bounded sets. 

Then 
 $\C_U = \{ f: \exists A\in U\, \exists k\, 
\forall \vec x\,\,  f(\vec x) \le h_A^{(k)}(\max\vec x)\} = 
\bigcup_{A\in U} \clmax{h_A}$.    

(Here, 
 $h_A^{(k)}$ is the $k$-fold iteration of the function $h_A$.)
\end{lemma}
\begin{proof} Write $\C_U':= \{ f: \exists A\in U\, \exists k\, 
\forall \vec x\,\,  f(\vec x) \le h_A^{(k)}(\max\vec x)\}$, 
$\C_U'' = 
\bigcup_{A\in U} \clmax{h_A}$.   

The inclusions $\C''_U \subseteq \C'_U \subseteq \C_U$ are trivial, and 
the inclusion  $\C_U'\subseteq \C_U''$ follows from the downward closure 
of $\clmax{h_A}$. 

To check $\C_U \subseteq \C_U'$,  it is 
enough to see that  $\C'_U$  is a clone.
 So let $f,g_1,\ldots, g_n\in \C'_U$, witnessed by~$A, A_1, \ldots, A_n, 
k, k_1, \ldots, k_n$.  Let $k^* = \max( k_1,\ldots, k_n)$, $A^* = 
A_1\cap \cdots A_n$. Then 
$$ f(g_1(\vec x), \ldots, g_n(\vec x)) \le h_{A}^{k} (\max(g_1(\vec
x), \ldots, g_n(\vec x))) \le h_{A}^{k}( h_{A^*}^{(k^*)}(\max\vec x))
\le h_{A\cap A^*}^{(k+k^*)}(\max\vec x)$$

\end{proof}

All clones that we construct in this paper will be of the form
$\C_U$ for some filter $U$. 

\np
\section{Creatures}

\subsection{Definitions}

\begin{definition}\label{def.planartree}
 A {\em planar tree} is a  tuple 
 $(T,{\trianglelefteq}, {<})$ where
\begin{itemize}
\itm A  $T$  is a nonempty set. (Elements of trees are often called ``nodes''.)
\itm B $\trianglelefteq$ is a partial order on~$T$ in which
  every set $\{\eta: \eta \trianglelefteq\nu\}$ is well-ordered by
  $\trianglelefteq$.   
\\ (We take $\trianglelefteq$ to be reflexive, and write $\vartriangleleft$ for the corresponding irreflexive relation)
\itm C $<$ is an irreflexive partial order on~$T$ such that any two $\eta\not=\nu$ in~$T$ are $<$-comparable iff they are $\trianglelefteq$-incomparable. ($x\le y$ means $x<y \vee x=y$.)
\itm D Whenever $\eta\trianglelefteq \eta'$ and $\nu \trianglelefteq \nu'$, then $\eta < \nu$ implies $\eta' < \nu'$. 
\end{itemize}
\end{definition}

\begin{example}
Let $T$ be a downward closed set of nonempty (possibly transfinite) 
sequences of ordinals. 
Then $T$ admits a natural tree order $\trianglelefteq$:
$\eta\trianglelefteq \nu$ iff $\eta $ is an initial segment of~$\nu$. 
We also have a natural partial order $<$, namely, the usual
lexicographic order of sequences of ordinals
(where sequences $\eta\vartriangleleft \nu$ are $<$-incomparable). Thus
$(T,\trianglelefteq,<)$
is a planar tree. 

It is easy to see that every planar tree in which the  relation $<$
is well-founded is isomorphic to a 
planar tree as described in this example.   None of our trees will
contain infinite 
chains, so they could be represented using sets of
finite (or even: strictly decreasing) sequences  of ordinals. 

For notational reasons, however,  we will use a completely different
way to represent trees.  The problem with the particular
implementation described above  is that we will have to ``glue'' old
trees together to obtain new trees, see~\ref{glue}; this means that 
the roots of the old trees will no longer be roots in the new tree. Since
we want to view the old trees as substrees of the new trees, it is 
not reasonable to demand that roots are always 
sequences of length~$1$. 
\end{example}

\begin{notation}\label{tree.notation}   
Let $(T,\trianglelefteq,<)$ be a planar tree. 
\begin{itemize}
\item[$*$]
We call 
${\trianglelefteq} = {\trianglelefteq}^T$ the ``tree order'', and 
${<} = {<}^T$ the ``lexicographic order'' of~$T$. 
\item[$*$]
  For $\eta\in T$ we write $\Succ_T(\eta)$ or sometimes $\Succ(\eta)$
  for the set of all direct successors of~$\eta$: 
$$ \Succ_T(\eta):= \{\,\nu\in T: \eta = \max\{\nu': \nu'
  \vartriangleleft \nu\}\,\}$$
\item[$*$]
  $\ext(T)$, the set of {\em external nodes} or {\em leaves} of~$T$ is
  the set of all $\eta$ with $\Succ_T(\eta)=\emptyset$. \\
  $\intr(T):= T \setminus \ext(T)$ is the set of {\em internal
  nodes}. 
\item[$*$]
  We let $\Wurzel(T)$ be the set of minimal elements of~$T$ (in the 
tree order $\trianglelefteq$).  If
  $\Wurzel(T)$ is a singleton, we call its unique element
  $\wurzel(T)$. 
\item[$*$]
  A {\em branch} is a maximal linearly ordered
 subset of~$T$  (in the sense of $\trianglelefteq$). \\
  The tree~$T$ is called ``well-founded'' iff $T$ has no infinite
  branches, or equivalently, no infinite linearly ordered subsets.
  Equivalently, if $(T,{\trianglerighteq})$ is well-founded in the
  usual sense. 
\\
   If $T$ is well-founded, then there is a natural bijection between
   external nodes and branches, given by~$\nu \mapsto \{\eta\in T: \eta
   \trianglelefteq\nu\}$.
\item[$*$]
 For any $\eta\in T$ we let $\yy T^\eta:= \{\nu: \eta \trianglelefteq \nu\}$; this is again a planar tree (with the inherited relations $\trianglelefteq$ and $<$). 
\\
More generally, if $H $ is a set of pairwise
$\trianglelefteq$-incomparable nodes of~$S$) (often $H \subseteq
\Wurzel(S)$), then we define 
$$\yy S^H:= \{\eta\in S: \exists \gamma \in H\,\, \gamma
\trianglelefteq \eta\} = \bigcup_{\gamma \in H} \yy S^\gamma .$$ This
is again a planar tree, and $\Wurzel(\yy S^H) = H$.  \\ If $H =
\{\gamma\in \Wurzel(S): \gamma_0 < \gamma\}$ for some $\gamma_0\in 
\Wurzel(S)$, then we write
 $\yyy S^{\gamma_0}<$ for $\yy S^H$.
\item[$*$]
  A {\em front} is a subset of~$T$ which meets each branch exactly once. 
(Equivalently, a front is a maximal subset of~$T$ set that is linearly ordered by~$<$.) 

  For example, $\leaf(T)$ is a front, and $\Wurzel(T)$ is also a front. 
  If $F \subseteq \intr(T)$ is a 
 front, then also $\bigcup_{\eta\in F} \Succ_T(\eta)$ is a front. 
 \\
Let $\eta\in \intr(T)$,  $F \subseteq \zz T^\eta$.   We say that 
$F$  is a ``front
above $\eta$''  iff $F$ meets every
branch of~$T$ containing $\eta$.
  Equivalently, $F$ is a front above $\eta$ if $F$ is a front in~$\zz T^\eta $.
\\
(For example, $\Succ_T(\eta)$ is a front above $\eta$.)
\item[$*$]
All trees $S $ that we consider will satisfy $\ext(S) \subseteq \lambda$, so it makes sense to define the following notation: 
\begin{itemize}
\item Let $S$ be a tree with $\ext(S)\subseteq \lambda$, and let $\eta\in S$. Then 
$\min_S[\eta]:= \min(\ext(\yy S^\eta))$.
\item Similarly  $\sup_S[\eta]:=\sup(\ext(\yy S^\eta))$.
\end{itemize}
\end{itemize}
\end{notation}

When $<$  and $\trianglelefteq$ are clear from the context we may just
call the tree ``$S$''; we  may later write $\trianglelefteq^S$,
$<^S$ for the respective relations.

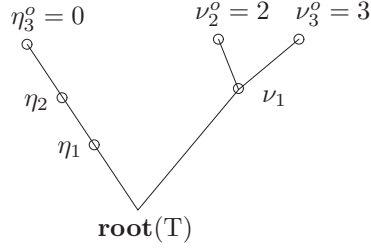
\begin{figure}
\setlength{\unitlength}{0.00035in}
\begingroup\makeatletter\ifx\SetFigFont\undefined%
\gdef\SetFigFont#1#2#3#4#5{%
  \reset@font\fontsize{#1}{#2pt}%
  \fontfamily{#3}\fontseries{#4}\fontshape{#5}%
  \selectfont}%
\fi\endgroup%
{\newcommand{\dashlinestretch}{30}
\begin{center}
\begin{picture}(5599,3614)(0,-10)
\put(530,3087){\ellipse{150}{150}}
\put(3387,3162){\ellipse{150}{150}}
\put(4587,3162){\ellipse{150}{150}}
\put(1053,2287){\ellipse{150}{150}}
\put(1535,1587){\ellipse{150}{150}}
\put(3687,2427){\ellipse{150}{150}}
\path(537,3087)(2187,612)(3687,2412)(3387,3162)
\path(3687,2412)(4587,3162)
\put(287,3287){\makebox(0,0)[lb]{$\eta_3^o=0$}}
\put(487,2087){\makebox(0,0)[lb]{$\eta_2$}}
\put(987,1387){\makebox(0,0)[lb]{$\eta_1$}}
\put(1587,127){\makebox(0,0)[lb]{\wurzel(T)}}
\put(4037,2162){\makebox(0,0)[lb]{$\nu_1$}}
\put(3037,3362){\makebox(0,0)[lb]{$\nu_2^o=2$}}
\put(4537,3362){\makebox(0,0)[lb]{$\nu_3^o=3$}}
\end{picture}
\end{center}
}
\caption{\label{examplefig} $\trianglelefteq$ and $\le$}
\end{figure}

We visualise such trees  as being embedded in the real plane ${\mathbb
R}^2$, with the order $\trianglelefteq$ pointing from the bottom to
the top, whereas the order $<$ can be viewed as pointing from left
to right.  (See Figure~\ref{examplefig},
 where we have $\eta_1\trianglelefteq\eta_2
\trianglelefteq\eta_3$, $\nu_1\trianglelefteq\nu_2 $,
$\nu_1\trianglelefteq\nu_3 $, $\nu_2 < \nu_3$, and $\eta_i < \nu_j$
for all $i,j\in \{1,2,3\}$.)

\begin{definition} \label{filter.converges}
 Let $(L,<)$ be a linear order, $D$ a filter on~$L$.  We say that
``$D$ converges to~$\sup L$'' iff for all $x_0\in L$ the set $\{y\in
L: x_0 < y\}$ is in~$D$. 
\end{definition}

\begin{fact} \label{no.last}
If  $(L,<)$ is a linear order, $D$ a filter on~$L$ converging to $\sup L$, then $L$ has no last element, and moreover, each $A\in D$ has no last element. 
\end{fact}
\begin{proof} If $A \in D$, $x_0\in A$, then the set $\{x\in A : x>x_0\}$ is 
in $D$ and hence cannot be empty. 
\end{proof}

\begin{definition}\label{def.abstract.creature}
 An abstract creature is a tuple $(S,{\trianglelefteq},{<}, D)$, where
\begin{itemize}
\itm A-D $(S,{\trianglelefteq},{<})$ is a planar well-founded tree 
\itm E $D = (D_\eta: \eta\in \intr(S))$ is a family of ultrafilters
\itm F For all $\eta\in\intr(S)$,
 the linear order $(\Succ_S(\eta), {<})$ has no
  last element.
\itm G For all $\eta\in\intr(S)$,
 $D_\eta$ is an ultrafilter on~$\Succ_S(\eta)$ which
  ``converges\label{converge}  to
  $\sup\Succ_S(\eta)$''
\end{itemize}
We sometimes write $(S,D)$ or just $S$ for creatures, if the other
parameters are clear from the context.   When an argument involves
several creatures $S, T, \ldots$, we may write $D^S$, $D^T$ etc for
the respective families of ultrafilters.   (The notation $D_S$ will 
be reserved for a quite different notion, see \ref{D.S}.)
\end{definition}

\begin{Remark}\label{induction}
 Since a creature $S$ is really a  well-founded tree
  $(S,{\trianglelefteq})$, we have that both 
  $(S,{\trianglelefteq})$ and   $(S,{\trianglerighteq})$ are
  well-founded. So when we prove theorems about the nodes of a
  creature $S$, or when we define a function on a creature, we can use one
  of two kinds of induction/recursion: 
\begin{itemize}
\item 
 ``Upward induction''.  Every nonempty $X \subseteq S$
  has a minimal element.  So if we want to define a function~$f$ 
``by recursion'' on~$S$ we may use the values of~$f\on \{\eta:
\eta\vartriangleleft \nu\}$ when we define $f(\nu)$.  Similarly, we can 
prove properties of all $\eta\in T$ indirectly by considering a
  minimal counterexample and deriving a contradiction.  
\item 
 ``Downward induction''. 
 Every nonempty $X \subseteq S$
  has a maximal element. So we can define a function $f$ on~$S$ by
 downward recursion ---  to define $f(\eta)$ we may use the function
 $f\on \{\nu: \eta\vartriangleleft \nu\}$, or more often the function
 $f\on \Succ(\eta)$. Similarly, we may use ``maximal counterexamples''
 in proofs of properties of all $\eta\in S$.  
\end{itemize}
\end{Remark}

\begin{motivation}
Mainly for notational reasons it will be convenient to be able to read
off information about the relations $\eta\trianglelefteq \nu$ and
$\eta < \nu$ directly from $\eta$ and $\nu$.  So we will restrict our
attention to a subclass of the class of all creatures:

 First we will
require all external nodes of our creatures to come from a fixed
linearly ordered set, the set of ordinals $< \lambda$.  We also 
require that the ``lexicographic'' order (see   \ref{tree.notation})
agrees with the  usual order of ordinals.

We then want to encode information about the location of any internal
node $\eta\in T$ within $T$ into the node $T$ itself.  It turns that
we can use the pair $( \min \ext\yy T^\eta, \sup \ext\yy T^\eta)$ as
``coordinates'' for~$\eta$.  Thus, all our creatures will be subsets
of $\lambda\cup \lambda\times\lambda$.
\end{motivation}

Definition \ref{Lambda} below is motivated by the following fact: 
\begin{Fact} Let $S$ be a creature with $\ext(S) \subseteq \lambda$. 
Then for all $\eta,\nu\in S$:\\
$\eta < \nu$ iff $\sup_S[\eta] \le \min_S[\nu]$.
\\
$\eta \vartriangleleft \nu$ iff $\min_S[\eta] \le \min_S[\nu]$ and 
$\sup_S[\nu] <  \sup_S[\eta]$. 
\end{Fact}
\begin{proof} 
 $\eta \vartriangleleft \nu$ implies 
that $\ext \yy T^\nu \subsetneq \ext \yy T^\eta$, so 
$ \min \ext\yy T^\eta\le  \min \ext\yy T^\nu$  and 
$\sup \ext\yy T^\eta \ge \sup \ext\yy T^\nu$.  In fact,
using \ref{def.abstract.creature}(d) it is easy to see that 
$\eta \vartriangleleft \nu$ even implies 
$\sup \ext\yy T^\eta >  \sup \ext\yy T^\nu$, so the map $\eta\mapsto
\sup[\eta]$ is 1-1. 
\end{proof}

\begin{defn}\label{Lambda}
Let $\Lambda:= \lambda \cup \{(i,j)\in \lambda\times\lambda:
i<j\}$.  
\\
We define two functions $\alpha$ and $\beta$ from $\Lambda$ into
$\lambda$: $\alpha(i,j)=i$, $\beta(i,j)=j$, $\alpha(i)=\beta(i)=i$ 
for all $i,j\in \lambda$.

We define two partial orders $<$ and $\trianglelefteq$ on $\Lambda$:
For all $\eta\not=\nu$:
$$\begin{array}{rclcl}
 \eta & < & \nu & \Leftrightarrow & \beta(\eta)\le \alpha(\nu)\\
 \eta & \vartriangleleft & \nu & \Leftrightarrow &\alpha(\eta)\le \alpha(\nu) \
 \wedge \  \beta(\eta)> \beta(\nu)\\
\end{array}
$$
\end{defn}

\begin{definition}\label{def.concrete}
 A concrete creature  (in the following just ``creature'') 
is a tuple $(S,{\trianglelefteq},{<}, D)$, where
\begin{itemize}
\itm A-G
      $(S,{\trianglelefteq},{<}^S, D^S)$ is an abstract creature
\itm H $S \subseteq \Lambda$, $<^S$ and $\vartriangleleft^S$
  agree with the relations   $<$ and $\vartriangleleft$ defined
  in~\ref{Lambda}.  
\itm I Each $\eta\in \intr(S)$ is a pair $\eta =
      (\alpha(\eta),\beta(\eta))$, and $\ext(S)\subseteq \lambda$. 
\itm J
  For all $\eta\in \intr(T)$: $\alpha(\eta)\le \min \ext(\yy
  T^\eta)$ and $ \sup\ext\yy T^\eta = \beta (\eta)$
\end{itemize}
\end{definition}

\begin{fact} 
Every creature (whose external nodes are a subset of~$\lambda $ with 
the natural order) 
is isomorphic to a concrete creature (replacing
each internal node  $\eta$ by the pair $(\min[\eta],\sup[\eta])$). 
\end{fact}

\begin{fact}
\label{global.order}
 If $S$ and $T$ are concrete creatures, and $\eta,\nu\in S\cap T$, then 
$\eta\trianglelefteq^S \nu$ iff $\eta\trianglelefteq^T \nu$,  and similarly
$\eta<^S \nu$ iff $\eta<^T \nu$. 
\end{fact}

We will often ``thin out'' creatures to get better behaved
subcreatures.  It will be easy to check that starting from a concrete
creature, each of the each of this thinning out processes will again
yield a concrete creature.

\np
\subsection{Small is beautiful}

\begin{definition}
Let $(S,D)= (S,\trianglelefteq,<,D)$ be a creature. 
We say that $(S,D)$ is 
\begin{ilist}
\item [\bf small,] if $\Wurzel(S)$ has a unique element: $\Wurzel(S) = 
  \{\wurzel(S)\}$. \\(This is a
  creature in the usual sense.)
\item [\bf medium,] if $\Wurzel(S)$ is infinite without last element 
but of cardinality $<
  \lambda$. \\ (Such a creature is often identified with the set (or
  naturally ordered sequence) $\{ \yy S^\gamma: \gamma \in
  \Wurzel(S)\}$ of small creatures.) 
\item [\bf large,] if $\Wurzel(S) \subseteq \lambda $ has size $\lambda$.  \\
(These creatures are usually called ``conditions'' in forcing
arguments.  They correspond to ``zoos'' in \cite{GoSh:808}. Again, it
may be convenient to identify such a large creature with a
$\lambda$-sequence of small creatures.) 
\end{ilist}
(We will not consider creatures $S$ with $1 < |\Wurzel(S)|$ where
$\Wurzel(S)$ has a last element.)
\end{definition}

\begin{fact}\label{nachfolger}
Let~$F$ be a front above $\eta$ (see~\ref{tree.notation}).  Assume $F \not=
\{\eta\}$. 
Then 
\begin{enumerate}
\item $F$ is linearly ordered by~$<$ and has no last element.
\item For all $\nu\in F$: $ \sup[\nu] < \sup [\eta]$
\item $\sup[\eta ] = \sup \{\sup[\nu]: \nu \in F\}$
\item $\sup[\eta ] = \sup \{\min[\nu]: \nu \in F\}$
\end{enumerate}
\end{fact}

\begin{minipage}[b]{9cm}
\begin{proof} 
We only show (1), the rest is clear.  Any two elements of~$F$ are
$\vartriangleleft$-incomparable, hence $<$-comparable.

Now let $\nu\in F$. 
We will find $\nu'\in F$, $\nu < \nu'$. 

Let $\eta \trianglelefteq \bar\eta \trianglelefteq \nu$, with
 $\bar\eta\in \Succ_T(\eta)$.  As~$\Succ_T(\eta)$ has no 
last element, we can find $\bar\eta'\in \Succ(\eta)$, $\bar \eta< \bar\eta'$. 

 So $\sup[\nu] \le \min[\bar \eta']$. 

There is~$\nu' \in F$ with $\bar\eta' \trianglelefteq \nu'$. By the definition of a planar tree, $\nu < \nu'$.  

\end{proof}
\end{minipage}
\qquad \begin{minipage}[b]{6cm}
\input etaprime.pstex_t 
\end{minipage}

 %
 %
 %
 %
 %

\subsection{Thinner creatures}

 \begin{dfact}[THIN] \label{thin}
 If $(S, D)$  is a small (or large)  creature, $S' \subseteq S$,
 then we write $S' \le\thin S$
 iff 
\begin{itemize}
\item $\Wurzel(S) \subseteq S'$.
\item  $\forall\eta\in S'\cap \intr(S): \succ_{S'}(\eta) \in D_\eta$.
\end{itemize}
In this case, $S'$ naturally defines again a small (or large, respectively)
creature $(S',D')$, by
letting   $D'_\eta := \{X \cap \succ_{S'}(\eta): X \in D_\eta\}$ for
all $\eta\in S'$, and by restricting $\trianglelefteq$ and $<$. 
 \end{dfact}

\begin{fact} If $S$ is a concrete creature, and $S'\le\thin S$, then
  also $S'$ is a concrete creature. 
\end{fact}
\begin{proof} Let $\eta = (\alpha,\beta)\in S'$.  We have to show that
  $\alpha \le \min_{S'}[\eta]$ and 
  $\beta =  \sup_{S'}[\eta]$.  The first property follows from 
$\alpha \le  \min_{S}[\eta]\le  \min_{S'}[\eta]$. 

For the second property we use downward induction.  Arriving at
$\eta$, we may assume  $\sup_{S'}[\nu] = \sup_{S}[\nu] $ for all 
$\nu\in \Succ_{S'}(\eta)$. Now 
 $\Succ_{S'}(\eta) $ is cofinal in 
$\Succ_S(\eta)$, hence also
 $\{\sup_{S'}[\nu]:\nu\in \Succ_{S'}(\eta) \} = 
\{\sup_{S}[\nu]:\nu\in \Succ_{S'}(\eta) \}$ is cofinal in 
$\{\sup_{S}[\nu]:\nu\in \Succ_{S}(\eta) \}$. 
\end{proof}

The following facts are easy: 
\begin{fact}\label{thinthin}
If $T $ and $S$ are small or large creatures, $T \le\thin S$, then for any
$\eta\in T$ we also have $\yy T^\eta \le\thin \yy S^\eta$. 
\end{fact}

\begin{fact} $\le\thin$ is transitive. 
\end{fact}

\subsection{drop, short, sum, glue}

\begin{dfact}[DROP]\label{drop}
Let $S$ and $T$ be large creatures.  We write $T\le\drop S$ iff
$\Wurzel(T)  \subseteq \Wurzel(S)$ (with the same order $<$)
and $ T = \yy S^{\Wurzel(T)}$.  (See~\ref{tree.notation}.)
\end{dfact}

Sometimes we drop only an initial part of the creature.  This relation
deservers a special name:

\begin{defn}[SHORT]\label{short}
Let $S$ and $T$ be large creatures.  We write $T\le\short S$ iff
there is some $\gamma \in \Wurzel(S)$ such that 
$\yyy S^\gamma< = T$. 

 We write $T\le\thinshort S$ iff there is some $T'$ with $T\le\short
 T'\le\thin S'$.   (Equivalently, if there is some $T'$ with $T\le\thin
 T'\le\short S'$.)
\end{defn}

\begin{definition}[SUM]
Let $(S,D)$ be a medium concrete creature. (See figure~\ref{sum}.) 
\begin{figure}
\input sum.pstex_t
\caption{A medium creature $S$\label{sum}}  
\input sum2.pstex_t
\caption{\label{sum2}A small creature  $\sum S$}  
\end{figure}

 Let $U$ be an ultrafilter on
$\Wurzel(S)$ converging to~$\sup\Wurzel(S)$
(see~\ref{filter.converges}).  Let 
$$ \alpha := \min \{ \alpha(\eta): \eta\in S\}, \quad 
  \beta := \sup \{ \beta(\eta): \eta\in S\} = \sup\ext(S), 
\quad \gamma:= (\alpha, \beta).$$ 
(Note that $\gamma \vartriangleleft \eta$ for all $\eta\in S$.)
 
Then $\sum (S,D ) = \sum\limits_{U} (S,D) =
\sum\limits_{U} S  $ is defined as the
following small concrete creature $(T,E)$  (see figure~\ref{sum2}): 
\begin{itemize}
\item[--] 
    $T:= \{\gamma\}\cup S$, $\wurzel(T)=\gamma $, $D_\gamma = U$, 

\item[--] 
 For all $\eta\in \Wurzel(S)=\Succ_T(\eta)$: $\yy T^\eta = \yy
 S^\eta$. 
\end{itemize}
\end{definition}

\begin{figure}
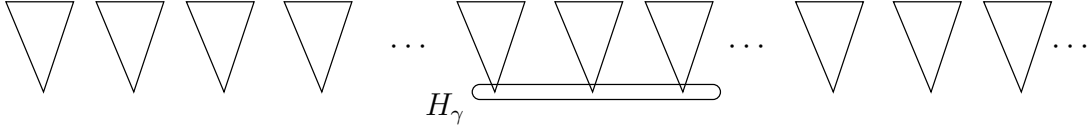

\input glue1.pstex_t 
\caption{\label{glue1}A large creature $S$}
\end{figure}

\begin{figure}
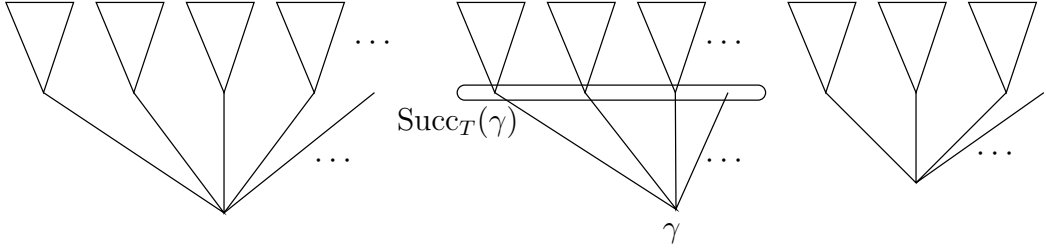
  
\input glue2.pstex_t 
\caption{\label{glue2}A large creature $T \le\glue S$} 
\end{figure}

\begin{definition}[GLUE]\label{glue}
Let $S$ and $T$ be  large concrete creatures.

We write 
 $T \le\glue S$ iff  for each $\gamma \in \Wurzel(T)$ the set
 $H_\gamma := \Succ_T(\gamma)$ is an interval in~$\Wurzel(S)$ with no
 last         
 element, and each $\yy T^\gamma $ can be written as $\sum_{U_\gamma,
 \gamma } \yy S^{H_\gamma}$ for some ultrafilters $U_\gamma$ (see
 figures~\ref{glue1} and~\ref{glue2}).

\end{definition}

\np

\subsection{Partition theorems}

\begin{lemma}\label{unary}
Let $E$ be a finite set. 
\begin{enumerate}
\item 
      If $S$ is a large or small creature, $c:S\to E$, then there is a
      creature $T \le\thin S$ such that $c\on \Succ_T(\eta)$  is constant
      for all $\eta\in T$. 
\item 
      If $S$ is a small creature, $c:\ext(S)\to E$, then there is a
      small creature $T \le\thin S$ such that $c\on \ext(T)$ is constant.
\item 
      If $S$ is a large  creature, $c:\ext(S)\to E$, then there are 
      large  creature $T' \le\drop\le T\le\thin S$
      such that $c\on \ext(T')$ is constant.
\end{enumerate}
\end{lemma}

\begin{proof}[Proof of (1)]
We define $T$ by upward induction, starting with $\Wurzel(T)=\Wurzel(S)$. 
Given $\eta\in T$, we find a set $A_\eta\subseteq \succ_S(\eta)$, 
$A_\eta\in D^S(\eta)$, such that $c\on A_\eta$ is constant, and we
let $\succ_T(\eta):=A_\eta$.
\end{proof}
\begin{proof}[Proof of (2)]
We define a map $\bar c:S\to E$ by downward induction (see \ref{induction}): 
\begin{enumerate}
\item[$*$]
     For $\eta\in \ext(S)$:    $\bar c(\eta)=c(\eta)$. 
\item[$*$] 
   For $\eta\in \intr(S)$ we find a (unique) value $e_\eta\in E$
   such that the set 
      $\{\nu\in \succ_S(\eta): \bar c(\nu)=e_\eta\}$ is in~$D_\eta$, and 
we set $\bar c(\eta):= e_\eta$. 
\end{enumerate}
Now we let $e_0:= \bar c(\wurzel(S))$, 
$$T:= \{\nu\in S: \forall \eta\trianglelefteq\nu\,  \bar c(\eta) = e_0\}.$$ 
Clearly $T \le\thin S$, and $c\on \ext(T)$ is constant with value $e_0$.
\end{proof}
\begin{proof}[Proof of (3)]

We apply (2) to each $\yy S^\gamma $, for all $\gamma\in \Wurzel(S)$
to get a large $T\le\thin S$ such that $c\on \ext(\yy T^\gamma)$ is
constant, say with value $e_\gamma$, for all $\gamma\in \Wurzel(T)$.
Now find $e_0$ such that the set $\{\gamma: e_\gamma = e_0\}$ has
cardinality $\lambda$, and let $T':= \bigcup_{e_\gamma=e_0} \yy
T^\gamma$.  Then $T'\le\drop T\le\thin S$, and $c$ is constant (with
value $e_0$) on $\ext(T')$.
\end{proof}

\np
\subsection{Comparing large creatures}

The constructions ``glue'', ``drop'' and ``thin'' are ways to get new,
in some sense ``stronger'' large creatures from old ones.  The
following definition gives a common generalization of the above
constructions.

\begin{definition}\label{le.creat}
Let $S,T$ be creatures. 
   We say $T \le   S$  iff there is a
 front~$F \subseteq T$   such that
\begin{itemize}
\item $F \subseteq \Wurzel(S)$
\item  for each    $\gamma \in F$:  
       $\yy T^\gamma  \le\thin \yy S^{\gamma }$ (see \ref{thin}). 
\end{itemize}
\end{definition}
\begin{figure}
\input le.pstex_t
\end{figure}
\begin{remark}
We usually consider this relation if both $S$ and $T$ are large, or
both are small, but we also allow the possibility that $S$ is large and $T$ is small. 

It is easy to see that if $S$ is small, 
and $T\le S$, then also $T$ must be  small and $T\le\thin
S$. 
\end{remark}

\begin{fact}\label{le.eta}
Assume that 
$T\le S$ are concrete creatures.  Then: 
\begin{enumerate}
\itm 1 For all $\eta\in T\cap S$ we have
 $\yy T^\eta \le\thin \yy S^\eta$. 
\itm 2 $\ext(T) \subseteq \ext(S)$, and $S$ is downward closed in $T$.
\end{enumerate}
\end{fact}

The next fact is the main reason for our notational device of
 ``concrete'' creatures (in~\ref{def.concrete}):
Thanks to \ref{global.order}, we may just  write 
 $\eta\trianglelefteq \nu$ in the proof below, rather than having to
 distinguish $\trianglelefteq^{S_1}$,
 $\trianglelefteq^{S_2}$, etc.

\begin{Fact}[Transitivity]
\label{le.transitive}
If $S_3\le S_2\le S_1$ are concrete creatures, then $S_3\le S_1$. 
\end{Fact}
\begin{proof}  Assume $S_3\le S_2\le S_1$,  where 
$S_k \le S_{k-1}$ is
  witnessed by a front  $F_k \subseteq S_k$ for $k=2,3$.

We claim that $F_2\cap S_3$ witnesses $S_3\le S_1$. Clearly $F_2\cap S_3
\subseteq \Wurzel(S_1)$.   To check that $F_2\cap S_3$
 is a front in~$S_3$, consider
any branch $b$ in $S_3$.   $b$ is of the form 
 $b=\{\eta\in S_3: \eta\trianglelefteq \nu_0\}$ for some 
$\nu_0\in \ext(S_3)$.  The set 
$\{\eta\in S_2: \eta\trianglelefteq \nu_0\}$ 
is also a branch in~$S_2$, so
it meets $F_2$ (hence $F_2\cap S_3$, by~\ref{le.eta}) in a singleton. 

For any $\eta \in F_2$, $\yy S_2^\eta \le \yy S_1^\eta $. Let $\gamma
\in F_3$, $\gamma \trianglelefteq \eta$.  Then we have 
$\yy S_3^\gamma \le\thin \yy S_2^\gamma $, so by \ref{thinthin} 
also 
$\yy S_3^\eta  \le\thin \yy S_2^\eta \le\thin S_1^\eta$. 
\end{proof}

\begin{examples}
\begin{enumerate}
\item 
For any $\gamma\in \Wurzel(S)$ we have $\yy S^\gamma \le S$.
\item 
 $S\le S$ is witnessed by the front $\Wurzel(S)$. 
\item 
Assume that $T\le\drop S$ or $T\le\thin S$. 
 Then  
again~$\Wurzel(T)$ witnesses $T\le S$. 
\item 
Assume that $T$ is obtained from $S$ as in GLUE (\ref{glue}).  Then
the front $\bigcup_{\gamma \in \Wurzel(T)} \Succ(\gamma)$  witnesses $T\le S$. 
\end{enumerate}
\end{examples}

\begin{lemma}\label{front.le.front}
  Let $S$ and $T$ be large 
concrete creatures, $T\le S$. 
 \begin{enumerate}
 \item $\ext(T) \subseteq \ext(S)$. 
 \item If  $F\subseteq S$ is a 
 front of~$S$, then $F\cap T$ is a front of~$T$.
 \end{enumerate}

\end{lemma}
\begin{proof}
(1) is clear. 

For (2), note that
nodes in $F\cap T$ are incomparable, because they were incomparable in
$S$, and $S$ and $T$ are concrete. 

Every external node of~$T$ is also an external node of~$S$, so every 
branch of~$T$ contains a branch of~$S$.   Hence every branch of~$T$
meets~$F$.  
\end{proof}

\np


\section{Creatures and functions}

\subsection{Weak and strong nodes}

In this section we will consider  functions $f:\lambda^k\to
\lambda$.    We will write tuples $(x_1,\ldots, x_k)\in \lambda^k$
as $\vec x$. For $\alpha\in \lambda$ we write $\vec x < \alpha$
iff we have $\max(x_1,\ldots, x_k) < \lambda$, similarly for 
$\vec x\le \lambda$.

However, the use of~$k$-ary functions is only a technicality;  the reader 
may want to consider only the case $k=1$, and then conclude the general results
either by analogy, or by assuming that all clones under consideration 
are determined by their unary fragments (this is true if all clones
contain a certain fixed 1-1 function $p:\lambda\times\lambda\to \lambda$).

\begin{definition}[Weak and strong nodes]\label{weak}
Let $f:\lambda^k \to \lambda $ be a $k$-ary function.

Let $(S,{\trianglelefteq},{<},D)$ be a creature,  $\eta\in S$.   
\begin{enumerate}
\item If $\eta\in \ext(S)$, then we say that  $\eta$ is $f$-weak.
\item $\eta\in\intr(S)$ is $f$-weak (in~$S$) iff 
 there is a $\vec y \in
  \lambda^k$, $\vec y \le \min[\eta]$, $f(\vec y) \ge \sup[\eta]$. \\
(Alternatively, we may say that $\eta$ is weaker than~$f$, or that 
$f$ is stronger than~$\eta$.)
\item $\eta\in \intr(S)$ is $f$-strong (in~$S$) iff for all $ \vec y \in
  \lambda^k$ with $\vec y < \sup[\eta]$ we have $f( \vec y) <
  \sup[\eta]$.  (Alternatively, we may say that $f$ is weaker
  than~$\eta$, or that $\eta$ is stronger than $f$.)
\end{enumerate}
We say that~$T$ is~$f$-strong iff each $\gamma\in  \Wurzel(T)$ is $f$-strong.
\end{definition}
\begin{figure}
\input weak.pstex_t
\end{figure}

\begin{remark}\label{weak.up}
If $\eta\trianglelefteq \nu$, and $\eta$ is $f$-weak, then also $\nu$
is weak.  So weakness  is inherited ``upwards''.  

Strength is in general not inherited downwards, but the following holds: 
\begin{quote}
If $F$ is a front above $\eta$, and all $\nu\in F$ are $f$-strong, then 
also $\eta$ is $F$-strong.
\end{quote}
\end{remark}

\begin{Fact}\label{strong.le.strong}  
Let $S$ and $T$ be concrete creatures. 
 Assume that $\eta\in S\cap T$ is $f$-strong (or $f$-weak) in
  $S$, and $T\le S$.  Then $\eta$ is again $f$-strong (or
  $f$-weak, respectively) in~$T$. 

  Similary, if $S$ is $f$-strong, then (by \ref{weak.up}) so is~$T$.
\end{Fact}
\begin{proof}
Assume $\eta$ is $f$-strong in~$S$, so $f(\vec y)<\sup_S[\eta]$ for all
$\vec y<\sup_S[\eta]$.  Since  $\sup_T[\eta]=\sup_S[\eta]$, so
$\eta$ will also be $f$-strong in~$T$. 

Assume $\eta$ is $f$-weak in~$S$, so $f(\vec y)\ge \sup_S[\eta] = \sup_T[\eta]$ 
for some $\vec y\le \min_S[\eta]$.  Clearly $\min_S[\eta]\le \min_T[\eta]$,
so also $\eta$ will also be $f$-weak in~$T$. 
\end{proof}

 \begin{Fact}\label{find.gauge.1}
Let $S$ be a  (large or small) creature, $f\in \O$.
\begin{enumerate}
\itm 1 There is $T \le \thinl S$  such
   that each $\eta\in T$ is either $f$-strong or $f$-weak. 
\itm 2 Moreover, there is $T$ as above, such that also for each internal 
$\eta\in \intr(T)$ either all $\nu\in \Succ_T(\eta)$ are $f$-strong, or 
                          all $\nu\in \Succ_T(\eta)$ are $f$-weak.
\end{enumerate}
 \end{Fact}
\begin{proof}  We define $T\le\thin S$ by upward induction, starting with
  $\Wurzel(T) := \Wurzel(S)$.  Now for each $\eta\in T$ we consider two
  cases: 
\begin{enumerate}
\item $\eta$ is $f$-strong (in $S$).  In this case we define
  $\Succ_T(\eta):=\Succ_S(\eta)$.  By \ref{strong.le.strong}, $\eta$ will
  also be $f$-strong in~$T$.
\item For some $\vec y < \sup_S [\eta]$ we have $f(\vec y) \ge \sup_S[\eta]$.\\
  Recall that (by~\ref{nachfolger}), $\sup_S[\eta] =
  \sup\{\min_S[\nu]: \nu\in \Succ_S(\eta)\}$.   So we can find 
$\nu_0 \in \Succ_S(\eta)$ with $\vec y \le\min_S[\nu_0]$.   
\\
Define $\Succ_T(\eta):= \{\nu\in \Succ_S(\eta): \nu_0 < \nu\}$.  (Note
  that this set is in~$D^T_\eta$.)\\
This ensures that $\eta$ will be $f$-weak in~$T$. 
\end{enumerate}
This completes the definition of~$T$, proving (1).

(2) now follows from (1) together with \ref{unary}(1).
\end{proof}

\np
\begin{Fact} \label{find.strong}
 Let $S$ be a large creature, $f:\lambda^k\to \lambda$.
Then there is 
$T \le S$ which is~$f$-strong. 
\end{Fact}

\begin{proof} 
Using the regularity of $\lambda$, we
 can find a continuous increasing sequence of ordinals $(\xi_i: i
< \lambda)$ with the following properties:
\begin{enumerate}
\item[--] For all $i<\lambda$, all $\vec x< \xi_i$: $f(\vec x) <
    \xi_i$.
\item[--] 
  For all $i<\lambda$, all $\gamma\in \Wurzel(S)$:
If $\min_S[\gamma] < \xi_i$, there there is $\gamma'>\gamma$ in
  $\Wurzel(S)$ with $\sup_S[\gamma]< \xi_i$. 
\item[--] 
  For all $i<\lambda$, the set $[\xi_i,\xi_{i+1})\cap \ext(S)$ is nonempty.
\end{enumerate}
These conditions will ensure that for all $i<\lambda $ the set 
  $$\Gamma_i:=\{ \gamma\in \Wurzel(S) : \xi_i\le 
   \min\ext(\yy S^\gamma)< \sup \ext(\yy S^\gamma)\le \xi_{i+1}\}$$
 is infinite with no last element.

 Now obtain $T$ from~$S$ by gluing together each set 
 $\{\zz S^\gamma: \gamma\in \Gamma_i\}$ 
(see \ref{glue}) for each $i<\lambda$. 
\end{proof}

\np

\subsection{Gauging functions with creatures} 

This section contains the crucial point of our construction: the close
correspondence between the relation $f \in \clmax{g}$ and the relation
$f<_S g$.

\begin{definition}
Let $S$ be a large 
 creature, $f:\lambda^k\to \lambda$, $F \subseteq S$
a  front.  We say that $F$ gauges~$f$ (in $S$) if
\begin{itemize}
\item For all $\eta\in F$: $\eta$ is $f$-strong. 
\item Whenever $\eta \vartriangleleft \nu$, $\eta\in F$, 
 then $\nu$ is $f$-weak.
\end{itemize}
We say that $S$ gauges $f$ if there is a front $F\subseteq S$ gauging $f$. 
\end{definition}

\begin{Fact}
Let $T\le S$ be large concrete creatures. 
If $S$ gauges $f$, then also $T$ gauges $f$. 
\end{Fact}
\begin{proof}
By \ref{front.le.front}, $F\cap T$ is a front in $T$. 
Let $F \subseteq S$ gauge $f$ (in~$S$),
  then $F\cap T$ still gauges $f$ (in $T$), witnessing that $T$ gauges $f$. 

(Use \ref{le.eta}.)
\end{proof}

\begin{fact} \label{find.gauge.2}
For every function $f\in \O$ and every large creature $S$ which is
$f$-strong
there is a
large creature $T\le\thin S$ which gauges $f$. 
\end{fact}
\begin{proof}

By \ref{find.gauge.1}, we can first find $T\le\thin S$ such that
all nodes in $T$ are $f$-strong or $f$-weak, and that all internal nodes 
have either only $f$-weak successors, or only $f$-strong successors. 

Now let $F$ be the set of all $\eta\in \intr(T)$ with the property
\begin{quote} $\eta$ is $f$-strong, but all $\nu\in \Succ(\eta)$ are
  $f$-weak. 
\end{quote}
Every  branch $b$ of $T$ contains an $f$-strong node (in $\Wurzel(T)$) and
an $f$-weak node (in $\ext(T)$)
 so  $b$ contains a highest strong node $\eta_b$; since $\eta_b$ has 
 a weak successors, all successors of $\eta_b$ are weak;
 hence $\{\eta_b\}= b\cap F$.
 Hence
 $F$ is a front, and clearly $F$ gauges $f$. 
\end{proof}

\begin{defn}\label{h.F}
Let $S$ be a creature, $F\subseteq S$ a front. We let 
$$ \llim_F:= \{ \sup\nolimits_S[\eta]: \eta\in F\}$$
and we write $h_F$ for the function $h_{\llim_F}$. 
\end{defn}
\begin{remark}
In the special case that $F=\ext(S)$, we have $\llim_F = F$, so our
 (new) definition of~$h_F$ agrees with our (old) definition
 in~\ref{h.A} of~ $h_{\ext(S)}$.  However, we will usually only
 consider fronts $F\subseteq \intr(S)$.
\end{remark}

\begin{remark}
If $F$ contains only internal nodes, then each point of 
$U_F$  is a limit point of~$\ext(S)$. We will see below
that  $h_F$  grows much faster than $h_{\ext(S)}$.  

 In an informal sense, $h_F$ is the smallest function that is still
 stronger than each $\eta\in F$.  The next lemma captures a part of
 that intuition.
\end{remark}

\begin{defn} For any $A\subseteq \lambda $ we write $f\le_A g$ 
 $f\in \clmax{h_A,  g}$.
\end{defn}

\begin{fact} The relation $\le_A$ is transitive. \end{fact}

\begin{lemma}\label{h.F.g}
Let $S$ be a large creature, $F\subseteq S$ a front.
Let $g$ be a function which is stronger than 
 each $\eta\in F$.   Then $h_{F}\le\es g $. 
\end{lemma}
\begin{proof} 
Let  $A:= \ext(S)$. 

For each $\eta\in F$ fix $\vec x_\eta$ such that $\max(\vec x_\eta)\le
\min[\eta]$, and $g(\vec x_\eta) \ge \sup[\eta]$. 

We will define a function $\vec y: \lambda\to \lambda^k$: 
\begin{quote}
           For each  $\alpha\in A$ we can find
           $\eta = \eta_\alpha\in F$ with
           $\eta_\alpha\trianglelefteq \alpha$.  
	   Let $\vec y(\alpha) = x_{\eta_\alpha}$.
	   \\
           For $\alpha\in \lambda\setminus A$ let $\vec y(\alpha) = \vec 0$. 
\end{quote}
Clearly $\vec y(\alpha) \le  \alpha$, so the function $y$ is in~$\C_\mx$.

For $\alpha\in \ext(S)$ we have 
$$ h_F(\alpha)  = \sup[\eta_\alpha] \le g(\vec y(\alpha))$$
and for $\alpha \notin \ext(S)$ we have $h_F(\alpha) = h_F(h_A(\alpha))$. 
In any case we have 
$$ h_F(\alpha) \le  g(\vec y( h_A(\alpha))),$$
therefore $h_F\in 
\clmax{ h_{A}, g}$.
\end{proof}

\begin{figure}
\input front.pstex_t
\end{figure}

\begin{lemma}\label{f.h.F}
Let $S$ be a large creature, $F\subseteq S$ a front.
Let $f$ be a function weaker than all $\eta\in F$.
  Then $f\le\es h_F$. 
\end{lemma}
\begin{proof}
For any $\vec x$, let $\eta\in F$ be minimal such that $\sup[\eta] > \vec x$. 
Then $h_F(\max \vec x) = \sup[\eta]$, but (as $\eta$ is $f$-strong), 
$f(\vec x) < \sup [\eta]$.   Hence $f(\vec x) < h_{F}(\max\vec x)$
for all $\vec x$, so $f\in \clmax {h_F}$. 
\end{proof}

\begin{lemma}\label{h.F.f}
Let $S$ be a large creature, $F\subseteq \intr(S)$ a front.
Let $f$ be a function which is weaker than each $\eta\in F$.

Then  $h_F\not\le\es f$. 
\end{lemma}
\begin{proof}
Pick any  $\eta\in F$, 
 and let $\xi:= \sup[\eta]$. 
Let 
$$ \D:= \{ c\in \O: \forall \vec x < \xi \, (c(\vec x) < \xi)\}.$$
Then $\D$ is a clone containing  $f$ (as $\eta$ is $f$-strong). 
As $\xi$ 
is a limit point of~$\ext(S)$, we also have
 $h\es(\vec x) < \xi$ for all $\vec x<\xi$, so $h\es\in \D$.

Hence $\clmax{h\es,f} \subseteq \D$, but $h_F\notin \D$. 
So $h_F\notin \clmax{h\es,f}$.
\end{proof}

\begin{Notation}
\label{F.nu.notation}
If $F $ is a front in $\ext(S)$, $\nu\in \ext(S)$, then we
write $F^\cdot(\nu)$ for the unique $\eta\in F$ with
$\eta\trianglelefteq \nu$. 
\end{Notation}

Recall from \ref{weak.up} that ``higher'' nodes (in the sense of
$\trianglelefteq$) are usualler weaker (in the sense of $f$-weakness) 
than lower nodes.  This apparent reversal of inequalities lies at the
heart of the next definition.

\begin{Definition}\label{compare}
Assume that $S$ is a large creature gauging $f$ and $g$, witnessed by
fronts $F$ and $G$.  We write 
\begin{enumerate}
\item[] $f<_S g$ iff: For all $\nu\in \ext(S)$, $F^\cdot(\nu)$ lies
strictly \emph{above} $G^\cdot (\nu)$: $G^\cdot (\nu) \vartriangleleft
F^\cdot(\nu)$. (See \ref{F.nu.notation}.)
\item[] $f\sim_S g$ iff: For all $\nu\in \ext(S)$, 
    $F^\cdot(\nu)=G^\cdot (\nu)$.
\end{enumerate}

We say that ``$S$ compares $f$ and $g$'' iff $S$ gauges $f$ and $g$, and
one of 
$$ f<_S g\qquad  f\sim_S g \qquad g <_S f$$ 
holds. 
\end{Definition}

\np
\begin{Fact}\label{fg.TS}    
 If $f<_S g$, and $T \le S$, then $f<_T g$.  

Similarly:  If $f\sim_S g$, and  $T \le S$, then $f\sim_T g$.  
\end{Fact}

The following lemma is the core of the whole proof.

\begin{lemma}\label{less.generates}
Let $S$ be a large creature gauging $f$ and $g$. 

If $f <_{S} g$, then $f\in \clmax{h\es,g}$, but 
$g\notin \clmax{h\es,f\}}$. In other words: 
\begin{quote}
If $ f<_S g$, then $ f \le\es g $, but $ g \not\le\es f$. 
\end{quote}
\end{lemma}

\begin{proof}
Let $F$ gauge $f$.  So every $\eta\in F$ is $f$-strong but $g$-weak. 

By \ref{f.h.F}, we have $f\le\es h_F$  and by \ref{h.F.g} $h_F\le\es g$. 
So $f\le\es g$. 

If we had $g\le\es f$, then (as $h_F\le\es g$, by  \ref{h.F.g}) we would get 
$h_F\le\es f$, contradicting \ref{h.F.f}. 
\end{proof}

Lemma \ref{less.generates} shows that if $S$ can ``see'' that $g$
grows faster than $f$, then together with $h\es$, $g$ dominates $f$, but not
conversely.  We can also read this as 
\begin{quote}
If $f<_S g$, then ``on the set $\ext(S)$'' $g$ dominates $f$ quite strongly.
\end{quote}
But can we always find a creature $S$ that can compare 
the different behaviors of $f$ and $g$? This is answered in the next lemma. 

\begin{lemma}
Let $f,g\in \O$, and let $S$ be a large creature. Then there is a large creature $T\le S$ which compares $f$ and $g$. (See \ref{compare})
$$ f <_T g \mbox{\quad or\quad } f \sim_T g \mbox{\quad or\quad } g <_T f .$$
\end{lemma} 
\begin{proof} By \ref{find.strong} we can find $S_1\le S$ which is
  $f$-strong, and by \ref{find.gauge.2} we can find $S_2\le S_1$  gauging
  $f$, witnessed by a front $F$.   Similarly we can find $S_3\le S_2$
  gauging $g$, witnessed by $G$.  $F\cap S_3$ still witnesses that
  $S_3$ gauges also $f$. 

To each external node $\nu$ of $S_3$ we assign one of three colors, 
depending on whether 
\begin{enumerate}
\item  $F^\cdot(\nu) = G^\cdot(\nu)$,
\item or  $F^\cdot(\nu) \vartriangleleft  G^\cdot(\nu)$, 
\item or  $F^\cdot(\nu) \vartriangleright  G^\cdot(\nu)$
\end{enumerate}
Using \ref{unary} we can find $T\le S_3$ such that all branches of $T$ 
get the same color.   

Now $T\le S$, and one of $f\sim_T g$, $f<_T g$, $g<_T f$ holds. 
\end{proof}

\begin{fact} \label{max.fact}
Assume $f\sim_S g$ or $f<_S g$.   Let $F$ and $G$ be the fronts 
gauging $f$ and $g$, respectively. 

Then: 
\begin{enumerate} 
\itm 1 Every $\eta\in S$ which is $g$-strong is also $f$-strong. 
\itm 2 For all $\eta\in S$:  $\eta$ is $g$-strong iff $\eta$ is
	$\max(f,g)$-strong.
\itm 3 $G$ gauges $g$
\itm 4 $\max(f,g)\sim_S g$.
\end{enumerate}
\end{fact}
\begin{proof}
(1): On every branch in $S$ the $g$-strong nodes are exactly the nodes which are $\trianglelefteq G$; these nodes are $\trianglelefteq F$, hence $f$-strong.

(2): Let $\eta$ be $g$-strong, so $\forall \vec x < \sup[\eta]$ we 
	have $g(\vec x) < \sup[\eta]$. As $\eta$ is also $f$-strong, we 
	also have
	$$\forall \vec x < \sup[\eta]: \max(f,g)(\vec x) < \sup[\eta].$$

(3) By (2).

(4) By (3).
\end{proof}

\np
\section{fuzzy creatures}

Ideally, we would like to construct a decreasing sequence
$(S_i:i<\lambda^+)$ of creatures such that the relations $\bigcup_i
{<}_{S_i}$ and $\bigcup_i {\sim}_{S_i}$ can be used for the
construction described in \ref{proof.outline}. However, the partial
order $\le$ on creatures is not even $\sigma$-closed, i.e., we can find a countable decreasing sequence with no lower bound. 

We will now slightly modify the relation $\le$ between large creatures
 to a relation $\le^*$ which has better closure properties but still keeps the 
important properties described in \ref{less.generates}.

\subsection{By any other name: $\approx\thin$, $\approx\short$, $\approx$}

\begin{dfact}\label{fact.intersection}
Assume that $S$, $S_1$, $S_2$ are concrete creatures, and: 
\begin{ilist}
\item[either:]  $S$  is  small, and both $S_1$ and $S_2$ are $\le \thin S$,
\item[or:]  $S$ is large, and both $S_1$ and $S_2$ are $\le \thinshort S$. 
\end{ilist}

We define a structure  $T=(T, {\trianglelefteq}^T, {<}^T, D^T) $ 
(which we also call $S_1\cap S_2$) as follows: 
\begin{enumerate}
\item  $\Wurzel(T) = \Wurzel(S_1)\cap \Wurzel(S_2)$, $T= S_1\cap S_2$
\item 
       ${<}^T = {<^{S_1}}\cap {<^{S_2}}$, 
\item 
       ${\vartriangleleft}^T = 
            {\vartriangleleft^{S_1}}\cap {\vartriangleleft^{S_2}}$,
\item   $D_\eta = D^{S_1}_\eta \cap D^{S_2}_\eta $ for all $\eta\in T$
\end{enumerate}
Then  $T$ is a creature, and $T\le\thin S_1$, $T\le\thin S_2$  (or 
$T\le\thinshort S_1, S_2$, respectively). 
\end{dfact}

\begin{proof}
We first check that $T$ is a planar tree.  Clearly $T$ is nonempty, as
 $T$ contains $\Wurzel(S)= \Wurzel(S_1)=\Wurzel(S_2)$. 
 Hence we have \ref{def.planartree}(A). 

 The orders
 $\trianglelefteq^{S_1}$ and $\trianglelefteq^{S_2}$ agree on $T$, as they both are restrictions of $\trianglelefteq^S$, and the same is true for $<^{S_1}$ 
and $<^{S_2}$.  This implies \ref{def.planartree}(B),(C),(D). 

We now check that $T$ is a creature. 
 For any $\eta\in T$ and any $A \subseteq \succ(\eta)$ we have 
$$ A \in D^T_\eta \Leftrightarrow A \in D^{S_1}_\eta\,\wedge \, 
 A \in D^{S_2}_\eta \Leftrightarrow A\in D^{S}_\eta$$
so $D^T_\eta$ is indeed an ultrafilter, i.e., \ref{def.abstract.creature}(E).

  Using \ref{no.last} we see \ref{def.abstract.creature}(F),(G). 

$T\le S_1, S_2$ is clear. 
\end{proof}

\begin{definition} \label{approxt}
Let $S$, $S'$ be small or large creatures.    We write $S \approx\thin  S' $ 
for $$ \exists T: \  T \le\thin S \mbox{ and } T \le\thinl S'. $$

Let $S$, $S'$ be large creatures.  We write $S \approx S'$ for 
 $$ \exists T: \  T \le\thinshort S \mbox{ and } T \le\thinshort S'. $$ 
\end{definition}

\np
\begin{fact} $\approx\thin$ and $\approx$ are  equivalence relations.
\end{fact}

\begin{proof} 
If $S,S',S'',T,T'$ are small (or large) creatures such that  $T$
witnesses $S \approx\thin S'$ and   $T'$ witnesses $S'\approx\thin S''$,
 then by~\ref{fact.intersection} we  see that $T'':=
  T\cap T'$ is again a small (or large) creature, and $T''$ witnesses 
  $S\approx\thin S''$. 

The proof for $\approx $ is similar.
\end{proof}

\begin{defn}[The relation $\le^*$]\label{lestar}
Let $T$ and $S$ be large  concrete creatures.  We say that $T \le^* S$
if there is~$T' $ with~$T \approx T' \le S$. 
\end{defn}

\begin{lemma}[Pullback lemma]\label{find.t.prime}
If $T_1 \le S_1 \approx S_0$ are large creatures, then there is a large
creature  $T_0$ such that $T_1 \approx T_0 \le S_0$:
$$
\begin{array}{ccc}
 & & S_0\\
 & & \approx \\
T_1 & \le & S_1
 \end{array}
\qquad \Longrightarrow\qquad 
\begin{array}{ccc}
T_0 &\le & S_0\\
\approx & & \approx \\
T_1 & \le & S_1 
 \end{array}
$$ 
\end{lemma}

\begin{proof} Let $F$ witness $T_1 \le S_1 $, and let $\gamma_0\in
F $ be so large that for all $\gamma\in F$ with
 $\gamma  > \gamma_0$ we have  $\yy S_1^\gamma \approx\thin \yy S_0^\gamma$.  

Let $F_0:= \{\gamma \in F:  \gamma >\gamma_0\}$, and  
define 
$$ T_0 =  
\bigcup_{\gamma \in F_0}  \{ \eta\in T_1: \eta\trianglelefteq \gamma \} 
\cup \yy S_0^\gamma $$
$T_0$ can be naturally equipped with a creature structure
$(\trianglelefteq^{T_0},<^{T_0},D^{T_0})$
such that $T_0 \approx T_1$; for defining $D^{T_0}$ we use the
fact that for all
 $\eta\in T_0$ with $\eta\vartriangleleft  \gamma\in F_0$
 the set $\Succ_{T_0} (\eta)$ is either equal to
$\Succ_{T_1}(\eta)$, 
or an end segment of this set, so in any case is in~$D_\eta^{T_0}$. 

Now  clearly $T_0
\le S_0$ is witnessed by $F_0$. 
\end{proof}

\begin{Corollary}
The relation $\le^*$ (between large creatures) is transitive. 
\end{Corollary}
\begin{proof} Let $S_3\le^* S_2 \le^* S_1$. 
We use our ``pullback lemma''~\ref{find.t.prime}: 
$$
\begin{array}{ccccc}
      &   &   S'   &   \le   &   R   \\
      &   &  \approx \\
   T' &\le&   S  \\
  \approx\\
   T
\end{array}
\qquad \Longrightarrow\qquad 
\begin{array}{ccccc}
   T''   & \le  &   S'   &   \le   &   R   \\
   \approx   &   &  \approx \\
   T' &\le&   S  \\
  \approx\\
   T
\end{array}
$$
and then appeal to the transitivity of~$\le $ and $\approx$. 
\end{proof}

\np
\subsection{Fusion}

\begin{Lemma}\label{fusion.delta}
Let $\delta < \lambda $ be a limit ordinal.  
Assume that $(S_i: i < \delta)$ is a sequence of large concrete 
creatures satisfying 
$i<j \Rightarrow S_j \le^* S_i$.  

Then there is a large creature $S_\delta $ such that: for all $i< \delta $:
$S_\delta \le^* S_i$. 
\end{Lemma}
A main idea in the proof is to divide $\lambda $ into $\lambda$
many pieces, each of length $\delta $: $\lambda = \bigcup_{\xi <
\lambda } [ \delta \cdot \xi, \delta \cdot \xi + \xi )$.
\begin{proof}

By elementary ordinal arithmetic, 
 for  each $\zeta < \lambda$ there is a unique pair
$(\xi, i)$ with $\xi < \lambda$, $i  < \delta$, and 
$\zeta =   \dai$.

Recall the definition of large creatures: each internal node $\eta$ is a pair $(\alpha(\eta), \beta(\eta))$, and $\ext(\yy S^\eta)$
is a subset of the interval $[\alpha(\eta), \beta(\eta))$, with
  supremum 
$\beta(\eta)$.

We choose (inductively) a sequence $r(\zeta)$ (for $\zeta  <
\lambda$)
of roots such that 
for all $\xi < \lambda$, all $i < \delta$: 
\begin{itemize}
\item  $\rdai\in \Wurzel(S_i)$
\item For all $\zeta'<\zeta$:  $r(\zeta')< r(\zeta)$. 
\\{} [If $\zeta'=\delta\cdot \xi'+i'$,
       $\zeta= \delta\cdot \xi+i$ with $i\not=i'$,
       then $r(\zeta')\in S_{i'}$ and $r(\zeta)\in S_i$ 
       come from different creatures, but they can still be compared: 
       \\
                               $r(\zeta')< r(\zeta)$ means
			      $\sup_{S_i'}[r(\zeta')] \le  
			      \min_{S_i}[r(\zeta)] $.]
\end{itemize}

Considering the matrix $(\yy S_i^{\rdai)}: i <
\delta, \xi < \lambda )$ of small creatures, we first note that 
$$ \uu T_0\ := \ \bigcup_{ \xi < \lambda  } \bigcup_{i<\delta } 
 S_i^{\rdai)}$$
is a large concrete 
creature.  (Whenever 
$  \delta\cdot\xi' + i' <   \dai$, 
and $\eta'\in S_{i'}^{r(\delta\cdot\xi' + i')}$,
$\eta\in S_i^{\rdai}$, then $\eta'<\eta$.) 

We also see that  $\uu T_0 \le^* S_0$, because for each
 $\xi<\lambda $ and each $i< \delta $ there is a small creature~$X$
with $$\yy \uu T_0^\rdai = \yy S_i^\rdai \approx\thin X \le S_0.
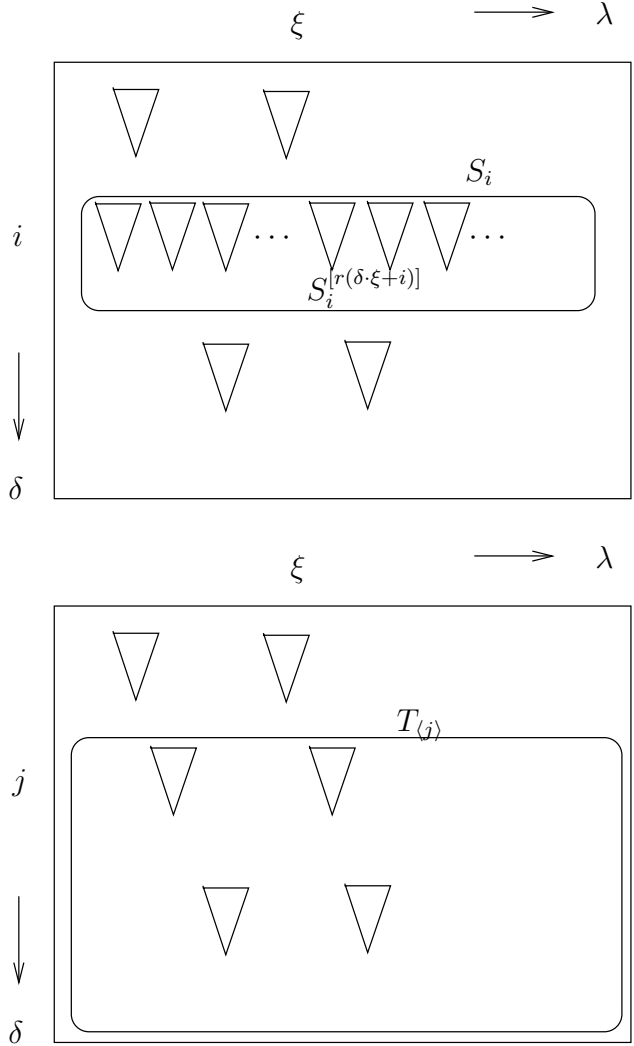
$$ 
\begin{figure}
\input matrix.pstex_t
\end{figure}

Similarly, we see that for every $j < \delta$ 
$$ \uu T_j\ := \ \bigcup_{ \xi < \lambda  } \bigcup_{j\le i<\delta } 
 S_i^{\rdai)}$$
is a large creature, and $\uu T_j \le^* S_j$. 
\begin{figure}
\input tj.pstex_t
\end{figure}

\newcommand{\tlim}{{\bar T}}

It remains to define a large creature $\tlim$ such 
that $\tlim \le^* \uu T_j$ for all $j<\delta $.

For each $\xi < \lambda $ the set 
$$ \uu T_{0,\xi}\ := \  \bigcup_{ i<\delta }  S_i^{\rdai)}$$
is a medium creature.

 Let $U_{0,\xi}$ be an ultrafilter on~$\Wurzel(\uu T_{0,\xi})$ which
 converges to $\sup\Wurzel(\uu T_{0,\xi})$, and let $r_\xi$ be a new root. 
Then $$\tlim_{0,\xi}:= \sum_{U_{0,\xi},r_\xi} \uu T_{0,\xi}$$
is a small creature, and $\tlim_0:= \bigcup_{\xi <
 \lambda}\tlim_{0,\xi}$  is a
 large creature.   By construction, $\tlim_0 \le\glue \uu T_0$.    

We can similarly define 
\begin{figure}
\input tjabar.pstex_t
\end{figure}
$$\uu  T_{j,\xi}\ := \  \bigcup_{j\le i<\delta }  S_i^{\rdai)}
\qquad  \tlim_{j,\xi}:= \sum_{U_{j,\xi},r_\xi} \uu T_{j,\xi}$$
(where  $U_{j,\xi}$ is the restriction of 
$U_{0,\xi} $ to $\Wurzel(\uu T_{j,\xi})$, an end
segment of~$\uu T_{0,\xi}$).

Again, $\tlim_j:= \bigcup_{\xi < \lambda} \tlim_{j,\xi}$  is a
 large creature satisfying $\tlim_j \le \uu T_j$.   

But by definition we have $\tlim_0 \approx\thin \tlim_j$, so $\tlim_0
\le^* \uu T_j$ for all $j< \delta $.

(See also figure~\ref{fusiond}.) 
\begin{figure}
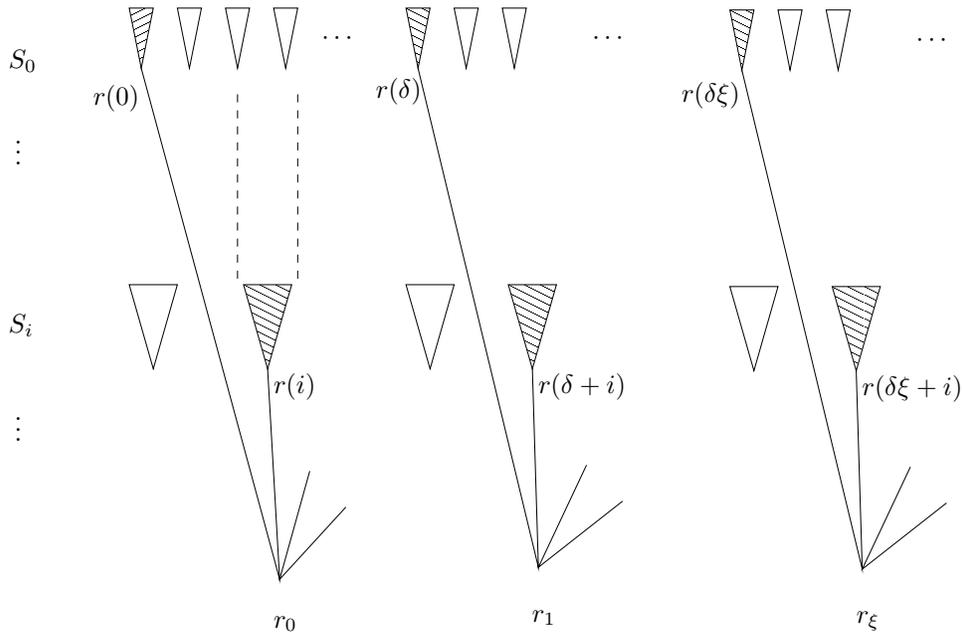
 
\input fusiond2.pstex_t 
\caption{\label{fusiond}A fusion of~$\delta $ many large creatures}
 \end{figure} 
\begin{figure}
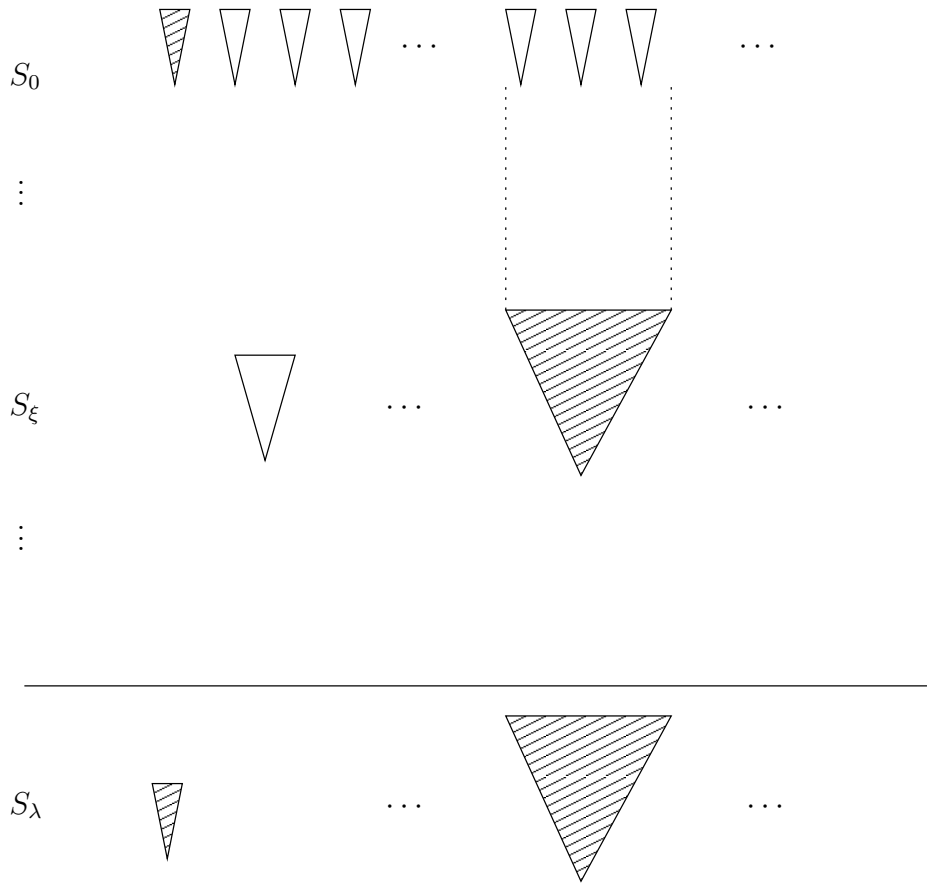

\input fusionl.pstex_t
\caption{\label{fusionl} A fusion of~$\lambda$ many large creatures}
\end{figure}

\end{proof}

\np

\begin{Lemma}\label{fusion.lambda}
Assume that $(S_\xi: \xi < \lambda)$ is a sequence of large 
concrete creatures satisfying 
$\xi<\xi' \Rightarrow S_{\xi '} \le^* S_\xi$.  

Then there is a large creature $S_\lambda  $ such that: for all $\xi< \lambda$:
$S_\lambda \le^* S_\xi$. 
\end{Lemma}
\begin{proof} 
We choose a fast enough increasing sequence
 $(r(\xi): \xi < \lambda)$
with $r(\xi)\in \Wurzel(S_\xi)$ such that: 
$$ \forall \zeta < \xi:  
r(\zeta) <  r(\xi).$$

Now let $T_0  := \bigcup_{\xi < \lambda }  \yy 
S_\xi^{r(\xi)}$, and similarly  
$T_\zeta  := \bigcup_{\zeta\le \xi < \lambda }  \yy
S_\xi^{r(\xi)}$. 
It is easy to see  $T_0 \approx T_\zeta \le S_\zeta$
for all~$\zeta$.  Hence $T_0\le^* S_\zeta$ for all $\zeta$.
\end{proof}

\begin{Corollary}\label{fusion.set}
Assume that ${\bf S}$ is a set of large concrete 
creatures which is linearly quasiordered by
$\le^*$, and assume that $|{\bf S}|\le \lambda$.  

Then there is a large creature $T$ with: $\forall S\in{\bf  S}: T\le^* S$. 
\end{Corollary}

\begin{proof} By \ref{fusion.delta} and \ref{fusion.lambda}. Use
  induction on~$|S|$. 
  \end{proof}

\np
\section{The filter $D_S$ and the clone $\C_S$} 

Recall that $f<_S g$ iff there are fronts $F,G\subseteq S$ gauging 
$f$ and $g$, respectively, such that $F$ meets each branch of $S$ below $G$.

\begin{Definition}
We write $f <_S^* g$ if there is $S'\approx S$, $f<_{S'} g$, similarly for $\sim^*$. 
\end{Definition}

\begin{Lemma}  If $f<_S^* g$, and $T\le^* S$, then $f <_T^* g$. 
\end{Lemma}
\begin{proof}
By the definition of $\le^*$ (see \ref{lestar}), there is $T_0$ such
that $ T \approx T_0 \le S$.  Let $S'\approx S$ be such that $S'$
gauges $f$.  Using the pullback lemma~\ref{find.t.prime}, we find
$T'\le S'$, $T'\approx T_0$.  So $T'\approx T$, $f<_{T'} g$ (by
\ref{fg.TS}), which implies $f<_{T}^* g$. 
$$
\begin{array}{ccc}
& & S'\\
& & \approx\\
T_0 &\le & S\\
\approx \\
T 
 \end{array}
\qquad \Longrightarrow\qquad 
\begin{array}{ccc}
T' &\le& S' \\
\approx  & & \approx\\
T_0 &\le & S\\
\approx & & \\
T & \\
 \end{array}
$$
\end{proof}

\np
\begin{defn}\label{D.S} Let $S$ be a large creature.   We define $$D_S:=\{ A
 \subseteq \lambda: \exists S' \approx S, \ext(S') \subseteq A\}.$$
\end{defn}

\begin{fact}\label{DS0} Let $S$ be a large creature. Then: 
 $A\in D_S$ iff there is $T\le\thinshort S$ with 
$\ext(T) \subseteq A$. 
\end{fact}
\begin{proof}
 If $A\in D_S$, then there are $S'$ and $T$ such that:
 \\$T\le\thinshort S$, \   $T\le\thinshort S'$, and 
 $\ext(S') \subseteq A$.  But then also 
$\ext(T)\subseteq \ext(S')\subseteq A$. 
\end{proof}

\begin{fact}\label{DS} Let $S$ be a large creature. Then: 
\begin{enumerate}
\itm 1 $D_S$ is a filter on $\lambda$, and all $A\in D_S$ are
unbounded.
\itm 2  If $S'\approx S$, then $D_S = D_{S'}$.
\itm 3 If $T\le S$, then $D_T \supseteq D_S$.
\itm 4 If $T\le^* S$, then $D_T \supseteq D_S$. 
\end{enumerate}
\end{fact}
\begin{proof}

(1) $D_S$ is clearly upward closed.  Let $A_1,A_2\in D_S$, witnessed
    by $S_1, S_2\le\thinshort S$, then $S_1\cap
    S_2$  witnesses $A_1\cap A_2\in D_S$.

(2) Immediate from the definition. 

(3) Follows from $\ext(T)\subseteq \ext(S)$ and the pullback lemma. 

(4) By (2) and (3). 
\end{proof}

\begin{defn} For any large creature $S$ we let 
$$\C_S:= \clmax{ h_A: A\in D_S} = \bigcup_{A \in D_S}\clmax{h_A}$$
\end{defn}
As a corollary to fact~\ref{DS} and fact~\ref{U.h.A} we get: 
\begin{fact}\label{CS} Let $S$ be a large creature. Then: 
\begin{enumerate}
\itm 1 $\C_S = \{f: \exists S'\approx S \, \exists k \, \forall \vec x\, 
(f(\vec x) \le h_{\ext(S')}^{(k)}(\max(\vec x))) \}$ 
\itm 2  If $S'\approx S$, then $\C_S = \C_{S'}$.
\itm 3 If $T\le S$, then $\C_T \supseteq \C_S$.
\itm 4 If $T\le^* S$, then $\C_T \supseteq \C_S$. 
\end{enumerate}
\end{fact}

\np
\begin{lemma} \label{58}
Let $S$ be a large creature, $f,g\in \O$, and assume $f<_{S}^* g$. 

Then $f\in \cl{\C_S \cup \{g\}}$, but 
 $g\notin \cl{\C_S \cup \{f\}}$. 
\end{lemma}
\begin{proof} 
  There is $S'\approx S$ with $f<_{S'} g$. But $D_S = D_{S'}$, so 
we may as well assume  $f<_S g$.

By \ref{less.generates}, $f\in \clmax{h\es, g}\subseteq 
 \clmax{ \{ h_A: A\in D_S\}\cup \{ g\}}  = \cl{\D_S\cup
 \{g\}}$. 

Assume that $g\in \cl{\C_S \cup \{f\}}$.  Then there is $A\in D_S$ 
such that $g\in \clmax{h_{A}, f}$. \\
 Let $S'\le\thinshort S$ with
$\ext(S') \subseteq A$.   Then 
$$g\in \clmax{h_{A}, f} \subseteq \clmax{h_{\ext(S')}, f}.$$
But  $S'\le S$ and   $f<_S g$ implies   $f<_{S'}g$, hence (again 
by \ref{less.generates}) we get  $g\notin \clmax{h_{\ext(S')}, f}$, 
a contradiction. 
\end{proof}

\section{Transfinite Induction}

\begin{defn} We say that a sequence $(S_i:i<\lambda^+)$ is
  ``sufficiently generic'' iff the sequence decreases with respect to
  $\le^*$: $$ \forall i < j:   S_j\le^* S_i,$$
and: 
$$ \forall f,g\in \O\, \exists i<\lambda^+:
     \ f<_{S_i} g \ \vee  f\sim_{S_i} g \
\vee  g<_{S_i} f$$ 
\end{defn}
\begin{Lemma} Assume $2^\lambda = \lambda^+$.  Then there is a
  sufficiently generic sequence. 
\end{Lemma}
\begin{proof} This is a straightforward transfinite induction:
There are 
$2^\lambda$ many pairs $(f,g)\in \O\times \O$.  By our assumption 
$2^\lambda=\lambda^+$ we can enumerate all these pairs as 
$$ \O \times \O = \{(f_i, g_i): i< \lambda^+\}.$$
We can now find a sequence $(S_i:i<\lambda^+)$ of large 
concrete creatures such that the following hold for all $i$:
\begin{itemize} 
\item If $i $ is a limit ordinal, then   $S_i\le^* S_j$ for all $j<i$.
\item  $S_{i+1}\le S_i$.
\item $S_{i+1}$  gauges $f_i$ and $g_i$
\item $S_{i+1}$ compares $f_i$ and $g_i$. 
              $g_i <_{S_{i+1}} f_i$.
\end{itemize}
\end{proof}

\begin{conclusion} Let $(S_i:i<\lambda^+)$ be a sufficiently generic
  sequence.  Define $\C_\infty:=\bigcup _i \C_{S_i}$.  This is an increasing
  union 
of clones, so also $\C_\infty$ is a clone. 

Let $f<_\infty g$ iff there is $i$ such that $f<_{S_i} g$, or equivalently, 
iff there is $i<\lambda^+$ such that $f<_{S_i}^* g$.  Define $f\sim_\infty g$
analogously. 

Then the properties (a)(b')(c') in \ref{proof.outline} are satisfied, so
\ref{proof.outline}(1)(2)(3) holds; moreover, for all $f\in \O$ there
is $g $  with  $f<g$, so $[\C,\O]$ has no cotaom. 
\end{conclusion}

\begin{proof}
(a) If $f<_\infty g$, then $f<_{S_i}^* g$ for some $i$. 
By \ref{58}, $f\in \cl{\C_{S_i}\cup\{g\}}$, so 
             $f\in \cl{\C\cup\{g\}}$.

(b') If $g\in \cl{  \C_\infty \cup \{f\}}$, then there is $i<\lambda^+$ such 
that $g\in \cl{  \C_{S_i} \cup \{f\}}$, as the sequence $(\C_{S_i})$ is 
increasing, by~\ref{CS}.

Choose $j>i$ so large that $S_j$ compares $f$ and $g$, so one of
$f<_{S_j} g$,
$f\sim_{S_j} g$,
$g<_{S_j} f$ holds. The first alternative is excluded by \ref{58}.

(c') follows from \ref{max.fact}.

Finally, let $f\in \O$.   Find $i<\lambda^+$ such that $S_i$ gauges~$f$.
Let $A:= \{ \sup_{S_i}[\gamma]: \gamma\in \Wurzel(S_i)\}$, and let $g:= h_A$. 
Then
\begin{itemize}
\item[$(*)$] Each $\gamma\in \Wurzel(S_i)$ is $f$-strong but $g$-weak.  
\end{itemize}

Now find $j>i$ such that $S_j$ compares $f$ and $g$. The possibilities 
$g <_{S_j} f$ and 
$f \sim_{S_j} g$ are excluded by $(*)$, so 
$f <_{S_j} g$, hence also $f<_\infty g$. 
\end{proof}

\newpage
\bibliography{other,listb,listx}

\bibliographystyle{lit-plain}
  
\end{document}